\documentstyle{amsppt}
\document

\topmatter
\title Perverse Sheaves on affine Grassmannians \\
 and Langlands Duality
\endtitle
\author Ivan Mirkovi\'c and Kari Vilonen
\endauthor
\address {Department of mathematics, University of Massachusetts, Amherst, MA
01002, USA}
\endaddress
\email mirkovic\@math.umass.edu
\endemail
\address{Department of mathematics, Brandeis University, Waltham, MA 02454,
USA}
\endaddress
\email vilonen\@math.brandeis.edu
\endemail
\thanks I. Mirkovi\'c was partially supported by NSF
\endgraf K.Vilonen was partially supported by NSA and NSF
\endthanks

\endtopmatter

\define\V{{\operatorname {Vec}_\Bbbk}}
\define\Ve{{\operatorname {Vec}_\Bbbk^\e}}
\define\oh{{\operatorname H}}

\define\ps{\operatorname {P}_{\cs}}
\define\p {\operatorname {P}}
\define\pg#1{\operatorname P_{\GO}(\G,#1)}
\define\GO{{G(\Cal O)}}
\define\GK{{G(\Cal K)}}
\define\G{{\Cal G}}


\define\bc{{\Bbb C}}
\define\e{{\epsilon}}
\define\cs{{\Cal S}}
\define\ct{{\Cal T}}
\define\cf{{\Cal F}}
\define\ca{{\Cal A}}

\define\ck{{\Cal K}}
\define\co{{\Cal O}}


%


\define\ra{\rightarrow}

\define\inj{\hookrightarrow}    

\define\bb{\underset}           
\redefine\aa{\overset}          

\define\barr{\overline}         
\predefine\ss{\S}               

\define\sub{\subseteq}          

\define\ch{\check}              
\define\tim{\times}             






\define\12{ \frac{1}{2} }       


\define\df{=_{\text{def}}}

\define\inv{{}^{-1}}

\define\Hom{\operatorname{Hom}}

\define\a1{{\A^1}}              


\redefine\AA{\Cal A}
\define\BB{\Cal B}

\define\GG{\Cal G}

\define\II{\Cal I}

\define\PP{\Cal P}

\define\VV{\Cal V}

\define\k{{\Bbbk }}

\redefine\A{\Bbb A }


\redefine\S{\Bbb S}

\define\Z{\Bbb Z}





\define\la{\lambda }




















\define\PGG#1{\operatorname P_{\GO}(\GG,#1)}








\define\htt{\operatorname{ht}}

\define\op{\operatorname{P}}

\define\lB{	\left(	}             	
\define\rB{	\right)	}

\subheading{\bf 1. Introduction}
\vskip .5cm

In this paper we outline a proof of a  geometric version of  the Satake
isomorphism. Namely, given a connected, complex algebraic reductive  group $G$
we show that the tensor category of representations of the dual group
$^L\!
G$ is naturally equivalent to a certain category of perverse sheaves on the
affine Grassmannian of $G$.
This can be extended to give
a topological realization of
 algebraic representations of
$^L\! G$
over any commutative ring $\k$ - the category of $\k$-representations of
$
^L\! G
$
is equivalent to a category of perverse
sheaves on the affine Grassmannian (a complex algebraic variety), with
coefficients in $\k$.

The above result over the complex numbers is not new. This case  has been
treated by Ginzburg in
\cite{Gi} and  some of the arguments in section 5 of  this paper are borrowed
from \cite{Gi}. However, at  crucial points our proof differs from Ginzburg's.
First, we use a more ``natural'' commutativity constraint for the
convolution product. This commutativity constraint, explained in section 3,
is due to Drinfeld and was explained to us by Beilinson. Secondly, in section
4,  we give a direct  proof that the global cohomology functor is
exact and  decompose this cohomology functor into a direct sum of weights
(Theorem 4.3).
The geometry underlying our arguments
leads to a construction of a canonical basis of Weyl
modules given by algebraic cycles.
Another consequence
is an explicit construction of the group algebra of $^L\! G$
in terms of the affine Grassmannian.
We completely avoid the use of the decomposition theorem of
\cite{BBD} which makes our techniques applicable to perverse sheaves with
coefficients over an arbitrary commutative ring $\k$.
In sections 1--5 we state the results when $\k$ is a field of characteristic
zero. The modifications needed for the general case are in
\ss 6. In
\ss 7
we derive the classical Satake isomorphism, for this we 
switch the setting to the
affine Grassmannian defined over a finite field
and   $\ell$-adic perverse sheaves.
This note contains indications of proofs of
some of the results. The details will appear elsewhere.

\vskip .5cm

\subheading{\bf 2. The Convolution Product}
\vskip .5cm

Let $G$ be a connected, complex algebraic reductive group. Denote by
$\Cal O = \Bbb C[[t]]$  the  ring of formal  power series in one variable and
by
$\Cal K = \Bbb C((t))$ its fraction field, the field of formal Laurent series.
The affine Grassmannian, as a set, is defined as $\G =
\GK/\GO$, where, as usual, $\GK$ and $\GO$ denote the sets of the
$\ck$-valued and the $\co$-valued points of $G$ respectively. The  sets
$\GK$,
$\GO$, and $\G$ have an algebraic structure as $\Bbb C$-spaces. The space
$\GO$ is a group scheme over $\Bbb C$ but  the spaces $\GK$ and $\G$ are only
ind-schemes\footnote{By an ind-scheme we mean an ind-scheme in a strict
sense, i.e., an inductive system of schemes where all maps are
closed embeddings.}. To see that
$\GK$ is an ind-scheme, one embeds $G$ in $SL_N(\bc)$. The filtration by order
of pole in $SL_N(\ck)$ induces a filtration of $\GK$ which exhibits
$\GK$ as an inductive limit of schemes. The filtration above is invariant under
the (right) action of $\GO$ on $\GK$ and thus, after taking the quotient of
$\GK$  by $\GO$ one gets a filtration of $\G$ which exhibits it as a union of
finite dimensional projective schemes. Furthermore, the morphism $\pi:\GK
\to
\G$ is locally trivial in the Zariski topology, i.e., there exists a Zariski
open subset $U\subset\G$ such that $\pi^{-1}(U)\cong U\times\GO$ and $\pi$
restricted to $ U\times\GO$ is simply projection to the first factor. For
details see for example \cite{BL1,LS}.

The group scheme $\GO$ acts on $\G$ with finite dimensional orbits. In order
to describe the orbit structure, let us fix a maximal torus $T\subset G$. We
write $W$ for the Weyl group and $X_*(T)$  for the coweights
$\operatorname{Hom}(\bc^*,T)$. Then the $\GO$-orbits on $\G$ are parameterized
by the $W$-orbits in $X_*(T)$, and given $\lambda\in X_*(T)$ the $\GO$-orbit
associated to it is $\G_\lambda = \GO
\cdot\lambda\subset\G$, where we have identified $X_*(T)$ as a subset of
$\GK$.

Let $\Bbbk$ be a field of characteristic zero, which we fix for the
sections 1--5.
All sheaves that we encounter in this paper will be sheaves in the
classical topology with the exception of \ss 7.  We denote by $\pg \Bbbk$
the category of
$\GO$-equivariant perverse
$\Bbbk$-sheaves  on $\G$ with finite dimensional support and by
$\ps(\G,\Bbbk)$ the category of perverse $\Bbbk$-sheaves on $\G$ which are
constructible with respect to the orbit stratification $\cs$ of $\G$ and which
have finite dimensional support. We use the notational conventions of
\cite{BBD} for perverse sheaves, in particular, in order for the constant sheaf
on a
$\GO$-orbit $\G_\lambda$ to be perverse it has to be placed in degree $-\dim
\G_\lambda$.

\proclaim{Proposition 2.1} The forgetful functor $\pg \Bbbk \to
\ps(\G,\Bbbk)$ is an equivalence of categories.
\endproclaim

We will now put a tensor category structure on $\pg \Bbbk$ via the convolution
product. Consider the following diagram of maps (of sets)
$$
\G\times\G @<p<< \GK \times \G @>q>> \GK \times_\GO
\G @>m>> \G\,.
\tag2.2
$$  Here $\GK \times_\GO\G$ denotes the quotient of $\GK \times \G$ by
$\GO$ where the action is given on the $\GK$-factor via right multiplication
by an inverse and on the $\G$-factor by left multiplication. The  $p$ and $q$
are projection maps and $m$ is the multiplication map. All other terms in (2.2)
have been given a structure of an ind-scheme except $\GK \times_\GO\G$. The
description of this structure is easier in the global context of section 3
where it is a special case of a more general construction and thus we
postpone the details. We define the convolution product
$A_1 *A_2$ of
$A_1,A_2\in\pg \Bbbk$ by the formula
$$ A_1 *A_2 \ = \ Rm_*\tilde A \qquad\text{where \ $q^*\tilde A  =
p^*(A_1\boxtimes A_2)$}\,.
\tag2.3
$$ To make sense of this definition we first use the fact  that $p$   and
$q$ are locally trivial in the Zariski topology. This guarantees the existence
of
$\tilde A\in\p_\GO(\GK\times_\GO\G,\Bbbk)$. To see the local triviality of
$q$ one can use the same arguments as for example in
\cite{BL1,LS}, and as was pointed out above, the local triviality of $p$ is
proved in those references. It remains to show that
$Rm_*\tilde A \in\pg\Bbbk$. To that end we introduce the notion of a
stratified semi-small map.

Let us consider two complex stratified spaces $(Y,\ct)$ and $(X,\cs)$ and a
map $f:Y\to X$. We assume that the two stratifications are locally trivial with
connected strata and that $f$ is a stratified with respect to the
stratifications $\ct$ and $\cs$, i.e., that for any $T\in\ct$ the image
$f(T)$ is a union of strata in $\cs$ and for any $S\in\cs$ the map
$f|f^{-1}(S):  f^{-1}(S)
\to S$ is locally trivial in the stratified sense. We say that
$f$ is a stratified semi-small map if

$$
\aligned a) \ \ &\text{for any $T\in\ct$ the map $f|\barr T$ is proper}
\\ b) \ \ &\text{for any $T\in\ct$ and  any $S\in\cs$ such that $S\subset
f(\bar T)$ we have}
\\ &\dim(f^{-1}(x)\cap \barr T) \leq \frac 1 2 (\dim f(\barr T) - \dim S)
\\ &\text{for any (and thus all) $x\in S$\,. }
\endaligned
\tag2.4
$$ Next the notion of a small stratified map. We say that $f$ is a small
stratified map if there exists a (non-trivial) open stratified subset $W$ of
$Y$ such that
$$
\aligned a) \ \ &\text{for any $T\in\ct$ the map $f|\barr T$ is proper}
\\ b) \ \ &\text{the map $f|W:W\to f(W)$ is proper and has finite fibers}
\\ c) \ \ &\text{for any $T\in\ct$, $T\subset W$, and  any $S\in\cs$ such that
$S\subset f(\barr T)-f(T)$}
\\ &\text{we have}\ \ \dim(f^{-1}(x)\cap \barr T) \leq \frac 1 2 (\dim
f(\barr T)
-
\dim S)
\\ &\text{for any (and thus all) $x\in S$\,. }
\endaligned
\tag2.5
$$

The result below follows directly from dimension counting:
\proclaim{Lemma 2.6} If $f$ is a semismall stratified map then
$Rf_*A\in\p_\cs(X,\Bbbk)$ for all $A\in\p_\ct(Y,\Bbbk)$\,. If $f$ is a small
stratified map then, with any $W$ as above, and any $A\in\p_\ct(W,\Bbbk)$, we
have $Rf_* j_{!*} A = \tilde j_{!*} f_*A$, where $j: W \hookrightarrow Y$ and
$\tilde j :f(W) \hookrightarrow X$ denote the two inclusions.

\endproclaim

We apply the above considerations, in the semismall case,  to our situation.
We take $Y=\GK\times_\GO\G$ and  choose $\ct$ to be the stratification whose
strata are $p^{-1}(\G_\lambda)\times_\GO\G_\mu$, for
$\lambda,\mu\in X_*(T)$\,. We also let $X=\G$, $\cs$ the stratification by
$\GO$-orbits, and choose
$f=m$. To conclude the construction of the convolution product on $\pg
\Bbbk$ it suffices to  note that the sheaf $\tilde A$ is constructible with
respect to the stratification $\ct$ and appeal to the following

\proclaim{Theorem 2.7} The multiplication map
$\GK \times_\GO \G @>m>> \G$ is a stratified semi-small map with respect to
the stratifications above.
\endproclaim

For an outline of proof, see remark 4.11.

One can define the convolution product of three sheaves completely analogously
to (2.3). This gives an associativity constraint for the convolution product
thus giving $\pg \Bbbk$ the structure of an associative tensor category. In
the next section we construct a commutativity constraint for the convolution
product.

\vskip 1cm

\subheading{\bf 3. The Commutativity Constraint}
\vskip .5cm

In order to construct the commutativity constraint we will need to consider
the convolution product in the global situation. Let $X$ be a smooth curve
over the complex numbers. Let $x\in X$ be a closed point and denote by
$\Cal O_x$ the completion of the local ring at $x$ and by $\Cal K_x$ its
fraction field. Then the Grassmannian $\G_x = G(\Cal K_x)/G(\Cal O_x)$
represents the following functor from $\Bbb C$-algebras to sets :
$$ R \mapsto \{
\cf
\text{  a  $G$-torsor on $X_R$,
$\nu: G\times X^*_R \ra \cf|X^*_R$ a trivialization on $X^*_R$  }
\}\,.
\tag3.1
$$ Here the pairs $(\cf,\nu)$ are to be taken up to isomorphism, $X_R =
X\times\text{Spec}(R)$, and $X^*_R= (X-\{x\})\times \text{Spec}(R)$\,. For
details see for example \cite{BL1,BL2,LS}. We now globalize this construction
and at the same time form the Grassmannian at several points on the curve.
Denote the
$n$ fold product by
$X^n = X\times
\dots\times X$  and consider the functor
$$ R \mapsto \left\{\aligned &(x_1,\dots,x_n) \in X^n(R), \ \ \cf\text{ a
$G$-torsor on
$X_R$\,, }\\ &\text{$\nu_{(x_1,\dots,x_n)}$ a trivialization of
$\cf$ on $X_R - \cup {x_i}$}\endaligned\right\}\,.
\tag3.2
$$ Here we think of the points $x_i: \text{Spec}(R) \to X$ as subschemes of
$X_R$ by taking their graphs. One sees that the functor in (3.2) is
represented by an ind-scheme $\G_X^{(n)}$. Of course  $\G_X^{(n)}$ is an
ind-scheme over
$X^n$ and its fiber over the point
$(x_1,\dots,x_n)$ is simply $\prod_{i=1}^k \G_{y_i}$\,, where
$\{y_1,\dots,y_k\}=\{x_1,\dots,x_n\}$, with all the $y_i$ distinct. We write
$\G_X^{(1)} = \G_X$.

We will now extend the diagram of maps (2.2), which was used to define the
convolution product, to the global situation, i.e., to a diagram of ind-schemes
over $X^2$:
$$
\G_X\times\G_X @<p<< \widetilde{\G_X\times \G_X} @>q>> \G_X  \tilde
\times
\G_X  @>m>>\G_X^{(2)}\,.
\tag3.3
$$ Roughly, the diagram starts with a pair of torsors, each trivialized off
one  point. One chooses a trivialization of the first torsor near the second
point, and uses it to glue the torsors.

More precisely,
$\widetilde{\G_X\times \G_X}$ denotes the ind-scheme representing the functor
$$ R \mapsto \left\{\aligned &(x_1,x_2)\in X^2(R); \ \cf_1,\cf_2\text{
$G$-torsors on
$X_R$;\
$\nu_i$ a trivialization of }
\\ &\text{$\cf_i$ on $X_R -x_i$, for $i=1,2$; \ $\mu_1$ a trivialization of
$\cf_1$ on
$\widehat{(X_R)}_{x_2}$}
\endaligned
\right\},
\tag3.4
$$ where $\widehat{(X_R)}_{x_2}$ denotes the formal neighborhood of
$x_2$ in
$X_R$. The \lq\lq twisted product" $\G_X  \tilde \times \G_X $ is the
ind-scheme representing the functor
$$ R \mapsto\left\{\aligned &(x_1,x_2)\in X^2(R); \ \cf_1,\cf \text{
$G$-torsors on $X_R$; $\nu_1$ a trivialization }
\\ &\text{of
$\cf_1$ on $X_R-x_1$;} \ \eta : \cf_1 | (X_R -x_2) @>{\ \ \simeq \ \ }>>
\cf| (X_R -x_2)
\endaligned
\right\}\,.
\tag3.5
$$ It remains to describe the morphisms $p$, $q$, and $m$ in (3.3). Because all
the spaces in (3.3) are ind-schemes over $X^2$, and all the functors  involve
the choice of the same $(x_1,x_2)\in X^2(R)$ we omit it in the formulas below.
The morphism $p$ simply forgets the choice of $\mu_1$, the morphism
$q$ is given by the natural transformation
$$
 (\cf_1,\nu_1,\mu_1;\cf_2,\nu_2) \mapsto (\cf_1,\nu_1,\cf,\eta),
\tag3.6
$$ where $\cf$ is the $G$-torsor gotten by gluing $\cf_1$ on $X_R - x_2$ and
$\cf_2$ on $\widehat{(X_R)}_{x_2}$ using the isomorphism induced by
$\nu_2\circ\mu_1^{-1}$ between $\cf_1$ and $\cf_2$ on
$(X_R-x_2)\cap \widehat{(X_R)}_{x_2}$. The morphism
$m$ is given by the natural transformation
$$ (\cf_1,\nu_1,\cf,\eta) \mapsto (\cf,\nu)\,,
\tag3.7
$$ where $\nu = (\eta \circ \nu_1)|(X_R - x_1-x_2)$.

Next, the global analog of $\GO$ is the group-scheme $G_X^{(n)}(\Cal O)$ which
represents the functor
$$ R \mapsto \left\{\aligned &(x_1,\dots,x_n) \in X^n(R), \ \ \cf\text{ the
trivial
$G$-torsor on
$X_R$\,, }\\ &\text{$\mu_{(x_1,\dots,x_n)}$ a trivialization of
$\cf$ on $\widehat{(X_R)}_{(x_1\cup\dots\cup x_n)}$}\endaligned\right\}\,.
\tag3.8
$$

Just  as in section 2 we define the convolution product of $\BB_1,\BB_2\in
\PP_{G_X(\Cal O)}(\G_X,\Bbbk)$ by the formula
$$
\BB_1 \bb{X}\to* \BB_2 \ = \ Rm_*\tilde \BB \qquad
\text{where \ $q^*\tilde \BB = p^*(\BB_1\boxtimes \BB_2)$}\,.
\tag3.9
$$ Precisely as in section 2, the sheaf $\tilde \BB$ exists because $q$ is
locally, even in the Zariski topology, a product. Furthermore, the map $m$ is a
stratified small map -- regardless of the stratification on $X$. To see this,
let us denote by $\Delta \subset X^2$ the diagonal and set $U = X^2-\Delta$.
Then we can take, in  definition (2.5), as $W$ the locus of points lying over
$U$. That $m$ is small now follows as
$m$ is an isomorphism over $U$ and over points of $\Delta$ the map $m$
coincides with its analogue in section 2 which is semi-small by theorem 2.7.

Let us now, for simplicity, choose $X=\Bbb A^1$. Then the choice of a global
coordinate on $\Bbb A^1$, trivializes $\G_X$ over $X$; let us write $\rho :\G_X
\to \G$ for the projection. Let us denote $\rho^0 =
\rho^*[1] : \PP_\GO(\GG,\Bbbk) \to \PP_{G_X(\Cal O)}(\G_X,\Bbbk)$\,.  By
restricting
$\G_X^{(2)}$ to the diagonal $\Delta \cong X$ and to $U$, and observing that
these restrictions are isomorphic to $\G_X$ and to
$(\G_X\times\G_X)|U$ respectively, we get the following diagram
$$
\CD
\G_X @>i>> \G_X^{(2)} @<j<< (\G_X\times\G_X)|U
\\ @VVV @VVV @VVV
\\ X @>>> X^2 @<<<\ \ \ U\ \ \ .
\endCD
\tag3.10
$$

\proclaim{Lemma 3.11} For $\AA_1,\AA_2\in\PP_\GO(\GG,\Bbbk)$ we have
$$
\aligned &\text{a)}\qquad
\rho^0 \AA_1   \bb{X}\to * \rho^0 \AA_2 \ \cong \ j_{!*}\left( (
\rho^0\AA_1\boxtimes\rho^0\AA_2)|U \right)
\\ &\text{b)}\qquad\rho^0(\AA_1*\AA_2) \ \cong \  i^0 ( \rho^0 \AA_1
\bb{X}\to * \rho^0\AA_2 ) \,.
\endaligned
$$
\endproclaim Part a) of the lemma follows from smallness of $m$ and lemma 2.6.

Lemma 3.11 gives us the following sequence of isomorphisms:
$$
\gathered
\rho^0(\AA_1*\AA_2) \cong i^0j_{!*} \left(
(\rho^0\AA_1\boxtimes\rho^0\AA_2)|U \right)
\\
\cong  i^*j_{!*}((\rho^0\AA_2\boxtimes\rho^0\AA_1)|U)
\cong\rho^0(\AA_2*\AA_1)\,.
\endgathered
\tag3.12
$$ Specializing this isomorphism to (any) point on the diagonal yields a
functorial isomorphism between $\AA_1*\AA_2$ and $\AA_2*\AA_1$. This
gives us a commutativity constraint making $\PP_\GO(\GG,\Bbbk)$ into a
tensor category.

\remark{Remark 3.13} The construction of the commutativity constraint can be
carried out in a more elegant way as follows. We first observe that the image
of the embedding $\rho^0 =\rho^*[1] : \PP_\GO(\GG,\Bbbk) \to
\PP_{G_X(\Cal O)}(\G_X,\Bbbk)$ consists precisely of objects in
$\PP_{G_X(\Cal O)}(\G_X,\Bbbk)$ which are \lq\lq constant" along $X$. This
subcategory of $\PP_{G_X(\Cal O)}(\G_X,\Bbbk)$ coincides with
$\PP_{\tilde G_X(\Cal O)}(\G_X,\Bbbk)$, where $\tilde G_X(\Cal O)$ denotes the
semi direct product of $G_X(\Cal O)$ and the groupoid which consists of pairs
of points $(x,y)\in X\times X$ together with an isomorphism between the formal
neighborhood of $x$ and the formal neighborhood of $y$. Now
$\rho^0 =\rho^*[1] : \PP_\GO(\GG,\Bbbk) \to\PP_{\tilde G_X(\Cal
O)}(\G_X,\Bbbk)$ is an equivalence whose inverse is $i^0 = i^*[-1]$, where
$i: \G_x \hookrightarrow \G_X$ is the inclusion. If $X$ is an arbitrary smooth
curve then the functor $i^0: \PP_{\tilde G_X(\Cal O)}(\G_X,\Bbbk) \to
\PP_\GO(\GG,\Bbbk)$ still has meaning and is an equivalence of categories. It
is clear that the convolution product (3.9) gives us a convolution product on
the category $\PP_{\tilde G_X(\Cal O)}(\G_X,\Bbbk)$. Thus, we can give the
construction of the commutativity constraint in terms of $\PP_{\tilde G_X(\Cal
O)}(\G_X,\Bbbk)$ and $i^0$ without specializing to $X= \Bbb A^1$ and choosing
a global coordinate.
\endremark
\vskip 1cm

\subheading{\bf 4. The Fiber Functor}
\vskip .5cm

Let $\Ve$ denote the category of finite dimensional $\Bbb Z/2\Bbb Z$ -\,
graded (super) vector spaces over
$\Bbbk$. Let us consider the global cohomology functor as $\Bbb H^*:
\pg\Bbbk \to \Ve$, where we only keep track of the parity of the grading on
global cohomology.  Then:

$$
\gathered
\text{The functor  $\Bbb H^*: \pg\Bbbk \to \Ve$ is a tensor functor
}
\\
\text{with respect to the commutativity constraint of section 3.}
\endgathered
\tag4.1
$$

Writing $r$ for the map $r:\G_X^{(2)} \to X^2$, this statement is an immediate
consequence of :
$$
\aligned &\text{a)}\ \ Rr_*(\rho^0(A_1)*_X\rho^0(A_2))|U  \ \ \text{is the
constant sheaf}\ \
\Bbb H^*(A_1)\otimes \Bbb H^*(A_2)\,.
\\ &\text{b)}\ \ Rr_*(\rho^0(A_1)*_X\rho^0(A_2))|\Delta = \rho^0(\Bbb
H^*(A_1*A_2))
\\ &\text{c)}\ \ Rr_*(\rho^0(A_1)*_X\rho^0(A_2)) \ \ \text{is a constant sheaf}
\endaligned
\tag4.2
$$
The claims a) and b) follow from lemma 3.11. It remains to note that, in the
notation of formula (3.9), the sheaf $R(r\circ m)_*\tilde B$ is constant; this
implies c)

Let $\V$ denote the category of finite dimensional  vector spaces over
$\Bbbk$. To make $\Bbb H^*: \pg\Bbbk \to \V$ into a tensor functor we
alter, following Beilinson and Drinfeld,  the commutativity constraint of \ss 3
slightly.
We consider the constraint from \ss 3  on the category
$\pg\Bbbk\otimes\Ve$ and restrict it
to a  subcategory that we
identify with
$\pg\Bbbk$.
Divide $\GG$ into unions of connected
components $\GG=\ \GG_+\cup\GG_-$ so that the dimension of  $\GO$-orbits
is even in $\GG_+$ and odd in $\GG_-$. This gives a $\Z_2$-grading
on the category   $\pg\Bbbk$ hence a new $\Z_2$-grading
on $\pg\Bbbk\otimes\Ve$. The subcategory of even objects
is identified with $\pg\Bbbk$ by forgetting the grading.

\proclaim{4.3 Proposition} The functor $\Bbb H^*: \pg\Bbbk \to \V$ is a tensor
functor with respect to the above commutativity constraint.

\endproclaim

We now come to the main technical result of this paper. In order to state it we
will fix some further notation. We choose a Borel subgroup $B\subset G$ which
contains the maximal torus $T$. This, of course, determines a choice of
positive roots. Let $N$ denote the unipotent radical of $B$. As usual, we
denote by
$\rho$ half the sum of positive roots of $G$. For any $\nu\in X_*(T)$ we write
$\htt(\nu)$ for the height of $\nu$ with respect to $\rho$. The
$N(\ck)$-orbits on
$\G$ are parameterized by $X_*(T)$; to each $\nu\in X_*(T) = \Hom(\Bbb C^*,T)$
we associate the $N(\ck)$-orbit $S_\nu \df\ N(\ck)\cdot\nu$. Note
that these orbits are neither of finite dimension nor of finite codimension.

\proclaim{Theorem 4.4} a) For all $\ca\in\pg\Bbbk$ we have
$$
\oh^k_c(S_\nu,\ca) = 0 \ \ \text{if} \ \ k\neq 2\htt(\nu)\,.
$$ In particular, the functors $\oh^{2\htt(\nu)}_c(S_\nu,\ \,): \pg\Bbbk \to
\V$ are exact.

b) We have a natural equivalence of functors
$$
\Bbb H^* \ \cong \ \bigoplus_{\nu\in X_*(T)} \ \oh^{2\htt(\nu)}_c(S_\nu,\ )\ :
\
\pg\Bbbk
\to
\V
$$
\endproclaim

This result immediately gives the following consequence:

\proclaim{Corollary 4.5} The global cohomology functor $\Bbb H^*: \pg\Bbbk
\to
\V$ is exact.
\endproclaim

Here is a brief outline of the proof of theorem 4.4. Let us consider unipotent
radical
$\barr N$ of the Borel $\barr B$ opposite to $B$. The $\barr N(\ck)$-orbits on
$\G$ are parameterized by $X_*(T)$: to each $\nu\in X_*(T)$ we associate the
orbit
$T_\nu = \barr N(\ck)\cdot \nu$\,. Recall that the $\GO$-orbits are
parameterized by $X_*(T)/W$. The orbit $S_\nu$  meets
$\G_\lambda$ iff $\nu\ \in\barr{\G_\la}\cap X_*(T)$, then
$$
\aligned &\text{a)} \ \ \dim(S_\nu\cap\G_\lambda) \ = \ \htt(\nu+\lambda)
\qquad
\text{if $\lambda$ is chosen dominant}
\\ &\text{b)} \ \ \dim(T_\nu\cap\G_\lambda) \ = \ -\htt(\nu+\lambda)
\qquad
\text{if $\lambda$ is chosen anti-dominant}
\\ &\text{c) \ \ the intersections in a) and b) are of pure dimension}\,.
\endaligned
\tag4.6
$$  In proving estimates  a) - c) we use the fact that the boundary $\partial
S_\nu$ is given by one equation in the closure $\barr S_\nu$.  From the
dimension estimates
(4.6a,b) above we conclude immediately that
$$
\aligned &\oh^k_c(S_\nu,\ca)\  =\ 0 \ \ \ \text{if} \ k>2\htt(\nu)
\\ &\oh^k_{T_\nu}(\G,\ca) \ =\ 0 \ \ \ \text{if} \ k<2\htt(\nu)\,.
\endaligned
\tag4.7
$$ Theorem 4.4 follows immediately from (4.7) and the following statement:
$$
\oh^k_c(S_\nu,\ca) \ = \oh^k_{T_\nu}(\G,\ca) \ \ \ \ \ \text{for all $k$}\,.
\tag4.8
$$  To see (4.8) we use the fact that $N(\ck)$-orbits and $\barr N(\ck)$-orbits
are in general position with respect to each other.

\remark{Remark 4.9} The decomposition of functors in theorem 4.4b is
independent of the choice of $N$. In the case of $N$ and its opposite
unipotent subgroup  $\barr N$  the corresponding decompositions are explicitly
related by $\oh^k_{S_\nu}(\G,\ca) \cong
\oh^k_{T_{w_0\cdot\nu}}(\G,\ca)$, where $w_0$ is the longest element in the
Weyl group. From this, and (4.8), we conclude that we could state theorem 4.4
replacing the functors $\oh^{2\htt(\nu)}_c(S_\nu,\ \,)$ by the equivalent set
of functors $\oh^{2\htt(\nu)}_{S_\nu}(\G,\ \,)$, where
$\oh^{2\htt(\nu)}_c(S_\nu,\ \,) \cong
\oh^{-2\htt(\nu)}_{S_{w_0\cdot\nu}}(\G,\
\,)$.
\endremark

\remark{Remark 4.10} The decomposition of $\G_\lambda$ into
$N(\ck)$-orbits and
$\barr N(\ck)$-orbits is an example of a perverse cell complex. Perverse cell
complexes are the analogues of CW-complexes for computing cohomology of
perverse sheaves instead of the ordinary cohomology. In the case at hand we
are in the situation analogous to the one for CW-complexes where  the
dimensions of all cells are of the same parity.
\endremark

\remark{Remark 4.11}  Theorem 2.7 follows from the estimate:
$$
\dim[	m\inv (\nu)\ \cap\ (p\inv(\GG_\la)  \bb{\GO}\to\tim \GG_\mu)  ]
\le
\htt(\la+\mu+\nu)
\tag4.12
$$
for  coweights $\la,\mu, \nu\in X_*(T)$ such that
$\la$ and $\mu$ are dominant and
$\nu\in\barr{\GG_{\la+\mu} }$.
This claim follows from the  estimate below which is a formal
consequence of (4.6):
$$
\text{
for a $T$-invariant subvariety $Y\sub\ \barr{\G_\la}$\,,
}\ \
\dim(Y)\le\ \bb{\nu\in \barr Y\cap X_*(T)}\to\max\ \htt(\la+\nu).
\tag4.13
$$

\endremark

\vskip 1cm

\subheading{\bf 5. The dual group}
\vskip .5cm

We will now apply Tannakian formalism as in \cite{DM} to $\pg\Bbbk$ and the
functor $\Bbb H^*$. In sections 2 and 3 we have given a tensor product
structure on the category $\pg\Bbbk$ via convolution and we have given
functorial associativity and commutativity constraints for this tensor
product.  To see that $\pg\Bbbk$ is a rigid tensor category, we still must
exhibit the identity object and construct duals. The identity object is given
by the sky scraper sheaf supported on the point
$1\cdot\GO\in\G$ whose stalk is
$\Bbbk$. The dual
$A\spcheck$ of a sheaf $A\in\pg \Bbbk$ is given as follows. Consider the
following sequence  of maps
$$
\G @<{\ \pi\ }<< \GK @>{\ i\ }>> \GK @>{\ \pi\ }>> \G\,,
\tag5.1
$$ where $i$ is the inversion on $\GK$, i.e., $i(g) = g^{-1}$. We define an
equivalence
$$
\iota:\pg \Bbbk \to \pg \Bbbk \text{ by} \ \iota(A) = \pi_* \tilde A \text{
where } i^*\tilde A = \pi^*A\,.
$$ Then the dual $A\spcheck$ is given by $A\spcheck = \iota (\Bbb D A)$, where
$\Bbb D$ denotes the Verdier dual.

In 4.3 we showed that $\Bbb H^*:\pg \Bbbk \to \V$ is a tensor functor.
Corollary 4.5 says that $\Bbb H^*$ is exact and from this it is not hard to
deduce that it is also faithful. Thus, we have verified that
$\pg\Bbbk$ together with $\Bbb H^*$ constitutes a neutral Tannakian category
and by \cite{DM, theorem 2.11} we conclude:

\proclaim{Proposition 5.1} There exists an affine group scheme $\check G$ such
that the tensor category $\pg\Bbbk$ is equivalent to the (tensor category) of
representations of $\check G$. This equivalence is given via the fiber functor
$\Bbb H^*$.
\endproclaim

This result, as well as the result below, can also be found in \cite{Gi}.
We claim:

\proclaim{Proposition 5.2} The affine group
scheme $\check G$ is isomorphic
to the Langlands dual of $G$.
\endproclaim

To see this, one may argue as follows. First of all, as in \cite{Gi}
one easily finds that
$\check G$ is noetherian and connected.
By theorem 4.4b) we conclude that the dual torus
$\check T$ of $T$ is contained in $\check G$ and then one shows, as in
\cite{Gi}, that the torus $\check T$ is maximal.
As one can express the root datum of a reductive group in terms of
its irreducible representations one concludes, following \cite{Gi},  that
the maximal reductive quotient of $\check G$ is the dual group of $G$.
It remains to eliminate the unipotent radical of $\ch G$,
but one can explicitly construct the group algebra of  $\check G$
and it turns out to be of the same size as the
the group algebra of the dual group.
have certain non-trivial self extensions of objects and this can easily be ruled
out (this argument is due to Soergel).

A few remarks are in order. Because $\check G$ is reductive, one concludes
immediately that $\pg\Bbbk$ is semisimple. One can also see directly that
$\pg\Bbbk\cong \ps(\G,\Bbbk)$ is semisimple, for example from \cite{Lu,
theorem 11c}.

Let us make the statements of propositions 5.1 and 5.2 more concrete. Let
$\lambda\in X_*(T)/W = X^*(\check T)/W$. To $\lambda$ we can associate an
irreducible representation $V_\lambda$ of the Langlands dual group
$\check G$ on one hand, and a $\GO$-orbit $\G_\lambda$, and thus an
irreducible perverse sheaf $\Cal V_\lambda= j_{!*} \Bbbk_\lambda[\dim
\G_\lambda]$, $j:\G_\lambda\hookrightarrow \G$, on the other. Under the
equivalence of proposition 5.1 the sheaf $\Cal V_\lambda$ and the
representation
$V_\lambda$ correspond to each other. Furthermore,  the representation space
of $V_\lambda$ gets identified with the global cohomology of $\Cal V_\lambda$,
i.e., $V_\lambda = \Bbb H^*(\G,\Cal V_\lambda)$. This interpretation gives a
canonical basis for $\Cal V_\lambda$ as follows. From theorem 4.4, the fact
that $j_{!*}\Bbbk_\lambda[\dim
\G_\lambda] = \break ^p j_{!}\Bbbk_\lambda[\dim \G_\lambda]$, and (4.6c) we
conclude:
$$
\gathered
\Bbb H^k(\G,\Cal V_\lambda) \ = \  \bigoplus\Sb\nu\in X_*(T)\\k =
2\htt(\nu)\endSb
\oh^{2\htt(\nu)}_c(S_\nu,\Cal V_\lambda) \ = \\ \bigoplus\Sb\nu\in X_*(T)\\k =
2\htt(\nu)\endSb\oh^{2\htt(\lambda+\nu)}_c(S_\nu\cap\G_\lambda,\Bbbk)
\ = \ \bigoplus\Sb\nu\in X_*(T)\\k = 2\htt(\nu)\endSb
\Bbbk[\operatorname{Irr}(S_\nu\cap\G_\lambda)]\,.
\endgathered
\tag5.3
$$ Here $\Bbbk[\operatorname{Irr}(S_\nu\cap\G_\lambda)]$ denotes the vector
space spanned by the irreducible components of
$S_\nu\cap\G_\lambda$\,. Thus we get
$$ V_\lambda \ = \ \Bbb H^*(\G,\Cal V_\lambda) \ = \ \bigoplus\Sb\nu\in X_*(T)
\endSb \Bbbk[\operatorname{Irr}(S_\nu\cap\G_\lambda)]\,.
\tag5.4
$$ Note that the results above imply that the cohomology group $\Bbb
H^*(\G,\Cal V_\lambda)$ is generated by algebraic cycles.

\vskip 1cm

\subheading{\bf 6. The case of rings and arbitrary fields}
\vskip .5cm

In this section we briefly indicate how the results in this paper can be
extended to
the case of rings. To this end, let $\k$ be a commutative, unital,
N\"otherian ring of
finite global dimension. We write $\pg \k$ for the category of
$\GO$-equivariant
perverse sheaves on $\G$ with coefficients in $\k$. Recall that this means,
in particular,
that the stalks of the perverse sheaves are finitely generated $\k$-modules.
All the
results in sections 2-4 go through in this context
except that in order that the convolution  of two perverse sheaves be
perverse one needs one of the sheaves to be in the full subcatgeory
 $\op^{proj}_\GO(\G,\k)$
consisting of sheaves
$\cf$
such that $\Bbb H^*(\G,\cf)$ is projective over $\k$.
The Tannakian
formalism yields

\proclaim{6.1 Proposition} There exists an affine group scheme $\check G$ over
$\operatorname{Spec}(\k)$ such that the tensor category
$\op^{proj}_\GO(\G,\k)$ is
equivalent to the tensor category of representations of
$\check G$ which are projective as $\k$-modules.
\endproclaim

Let us write $^L \!G$ for the unique split reductive group scheme over
$\operatorname{Spec}(\Bbb Z)$ corresponding to the root datum dual to that
of $G$
and let $^L \!G_\k$ denote the corresponding group scheme over
$\operatorname{Spec}(\k)$. Then

\proclaim{6.2 Theorem} The group scheme $\check G$ is isomorphic to
$^L \!G_\k$.
\endproclaim

The main new point is to show that
$\check G$  is
reduced, this is not
a priori clear  when  $\k$ is not a field of
characteristic zero.
We analyze explicitly the pro-projective object which
represents the
functor $\Bbb H^*: \pg \k \to \{\k-\text{modules}\}$. Its Verdier dual
corresponds to the group algebra  of $^L\!G_\k$.
With the help of the projective one finds that
the equivalence (6.1)
canonically extends to an equivalence of abelian categories
$\op_\GO(\G,\k)$ and  the  category of representations of
$\check G$.
The proof that $\ch G$ is integral (hence reduced),
is based on the perverse-sheaf analogue of the fact that
the tensor product of Weyl modules has a filtration such that the
succesive subquotients are Weyl modules.

Each $\GO$-orbit $\ \GG_\la \aa{j}\to\inj\GG$, defines two
standard 
perverse sheaves 
$
\II^!_\la(\k)
\df\ ^p\oh^0\lB Rj_!\k_{\GG_\la}[\dim(\GG_\la)]\rB
$
and
$\
\II^*_\la(\k)
\df\ ^p\oh^0\lB Rj_*\k_{\GG_\la}[\dim(\GG_\la)]\rB
$.
The usual (middle perversity) intersection homology sheaf
$
\VV_\la(\k)
\df\ j_{!*}\k_{\GG_\la}[\dim(\GG_\la)]
$,
is the image of the canonical map
$
\II^!_\la(\k)\ra \II^*_\la(\k)
$. (All three coincide if $\k$ is a field of characteristic zero.)
These standard perverse sheaves 
are topological realizations of the Weyl module
$\Bbb H^*\lB\GG,\II^!_\la(\k)\rB$ and its dual 
$\Bbb H^*\lB\GG,\II^*_\la(\k)\rB=$ 
sections of a line bundle on the flag variety.
The arguments  in \ss 5  now yield  canonical bases of both
the Weyl module and its dual. (When $\k$ is a field of characteristic zero
these are two  bases of
$\Bbb H^*\lB\GG,\VV_\la(\k)\rB$.) 
We
conclude with the following

\proclaim{6.3 Conjecture} The stalks of $\II^!_\la(\Z)$ and the costalks 
of $\II^*_\la(\Z)$ have no torsion.
\endproclaim

\vskip 1cm

\subheading{\bf 7. The Satake isomorphism}
\vskip .5cm

In this section we briefly discuss the case of the classical Satake
isomorphism. Let $G$ be a split reductive group over
$
{\Bbb F_q}$.  The
construction of
the affine Grassmannian can be performed in this context exactly the same
way as over
$\Bbb C$, with
$\Cal O=
{\Bbb F_q}[[t]]$ and $\Cal K=
{\Bbb
F_q}((t))$.
We consider the category $\op_{G(\Cal O)}
(\Cal G,\barr{	\Bbb	Q_\ell	})$
of $\GO$-equivariant {\it pure} perverse
$\barr{	\Bbb Q_\ell	}$-sheaves  of weight zero  on $\G$. The convolution
product preserves this category and the results in
sections 2-5 go through in this context, resulting in the following:

\proclaim{7.1 Theorem} There is an equivalence of tensor categories
$$
\op_{G(\Cal O)}(\Cal G,\barr	{\Bbb	Q_\ell}	) \cong Rep(
{^L}\! G_{	\barr{	\Bbb Q_\ell	} 		}	).
$$
\endproclaim

Let us now pass to the Grothendieck groups on both sides of 7.1. On the
right hand side we get the representation ring
$\text{Rep}[	{^L}\!G_{ \barr{	\Bbb Q_\ell  }	}	]$.  It
remains to interpret the
Grothendick group of the  left hand side. Because we are considering the
category of pure sheaves only, passing from sheaves to functions  on
$\GG(\Bbb F_q)$ via the trace of Frobenius results in an isomophism of the
Groethendieck group of  $\op_{G(\Cal O)}(\Cal G,\barr	{\Bbb	Q_\ell}	)$  and
the
 spherical Hecke algebra $\Cal H $ of compactly supported
$G(\Bbb F_q[[t]])$-bi-invariant functions on
$G(\Bbb F_q((t)))$.
Thus theorem 7.1 yields an
isomorphism of
$\barr{\Bbb Q_\ell}$-algebras
$$
\Cal H \ \cong \ \text{Rep}[  {^L}\!G_{	\barr{\Bbb Q_\ell}	}	]\,.
\tag7.2
$$
This is the classical Satake isomorphism.

\remark{7.3 Remark} Some of the arguments for sections 2-5 have to modified
slightly
when working in the context of $\ell$-adic sheaves. This has been done, in a
way different from ours, earlier by Ngo in \cite{N}.
\endremark

\enddocument